\documentclass[a4paper,11pt]{article}
\usepackage{amssymb,palatino}
\newtheorem{theorem}{Theorem}

\newtheorem{lemma}{Lemma}

\newtheorem{proposition}{Proposition}
\newcommand{\sty}{\displaystyle}
\newcommand{\QED}{\begin{flushright} $\Box$ \end{flushright}}
\newcommand{\binomiale}[2]{\left( \begin{array}{c} #1\\#2 \end{array} \right)}

\begin{document}

\title{Improvements on the Cantor-Zassenhaus \\ Factorization Algorithm}

\author{Michele Elia\thanks{Politecnico di Torino, Italy}, ~~
 Davide Schipani   \thanks{University of Zurich, Switzerland}
}

\maketitle

\thispagestyle{empty}

\begin{abstract}
\noindent
After revisiting the Cantor-Zassenhaus polynomial factorization algorithm, we describe a new simplified version of it, which entails a lower computational cost. Moreover, we show that it can be used to
 find a factor of a fully splitting polynomial of degree $t$ over $\mathbb F_{2^m}$ with $O(\frac{2^m}{3^{t}})$ attempts and over $\mathbb F_{p^m}$ for odd $p$ with $O(\frac{p^m}{2^{t}})$ attempts.  
\end{abstract}


\vspace{2mm}

\vspace{8mm}

\section{Introduction}
The Cantor-Zassenhaus polynomial factorization algorithm (\cite{cantor}) is an efficient (polynomial-time)
 probabilistic algorithm for factoring polynomials over a finite field $\mathbb F_{p^m}$, that are the
 product of irreducible polynomials with a common degree $s$ and multiplicity one. When the multiplicity
 is above $1$, the factors can be separated by computing the greatest common divisor of the given
 polynomial and its formal derivative. If the irreducible polynomials have different degrees, the factors
 are separated by computing the greatest common divisors with polynomials of the form $x^{p^{mr}-1}-1$,
 starting from $r=1$, so as to obtain the product of all irreducible factors of degree $r=1,2,\ldots$
 (see e.g. \cite{jurgen}).
Thus standard methods can be used to reduce the problem to the above case. 


We will now 
introduce the Cantor-Zassenhaus factorization algorithm, 
providing a non-standard explanation which will be the basis for the rest of the paper: in the Sections below we will show how it can be improved, giving a new
 description with a more favorable estimate of its complexity and success rate. In fact this description leads us to consider a deterministic version of the algorithm, so that we will be concerned with the problem of establishing how many attempts are needed in the worst case to obtain a factor (with probability $1$) and what is the least degree of the polynomial such that a factor is found with at most a fixed number of attempts. 

Let $\sigma(z)$ be a polynomial of degree $t$ over $\mathbb F_{p^m}$ which is  
 a product of irreducible polynomials of degree $s$, i.e. $t=s \cdot d$. 

Let us assume that $s=1$ as a first instance and suppose that the trivial factor $z$ does not divide $\sigma(z)$.

We first deal with the case $p=2$, and following \cite{cantor} we assume that $m$ is even,
 otherwise we would consider a quadratic extension solely for the computations. 
If $\alpha$ is a known primitive element of $\mathbb F_{2^m}$, we define $\ell_{m} = \frac{2^{m}-1}{3}$
 and $\rho= \alpha^{\ell_m}$, which is thus a primitive cubic root of unity in the field 
 $\mathbb F_{2^m}$.

Let $c(z)$ be a non-constant polynomial over $\mathbb F_{2^m}$ of degree less than $t$, and let
 $$    a(z) = c(z)^{\ell_{m}}  \bmod \sigma(z)  $$
which is again a polynomial of degree at most $t-1$. Furthermore, we have
$$  (c(z)^{\ell_{m}}+1)(c(z)^{\ell_{m}}+\rho) (c(z)^{\ell_{m}}+\rho^2) =  c(z)^{2^{m}-1}-1 ~~. $$

Now, either $\gcd(c(z),\sigma(z))$ is non-trivial (and thus we already have a factor of $\sigma(z)$)
 or else $c(z)^{2^{m}-1}-1 = 0 \bmod \sigma(z)$. In this latter case, if we write
  $c(z)^{2^{m}-1}-1= Q(z)\sigma(z)+R(z)$ and specialize it in the roots $\{z_i\}$ of $\sigma(z)$, we see that $R(z)$, which is a polynomial of degree $t-1$, takes the value $0$ for all $t$ roots, as $\beta^{2^{m}-1}-1=0$ for any $\beta\in \mathbb F_{2^m}^*$. This implies that $R(z)$ is identically $0$.
Thus we can write
$$ (c(z)^{\ell_m}+1)(c(z)^{\ell_m}+\rho) (c(z)^{\ell_m}+\rho^2) = (a(z)+1)(a(z)+\rho)(a(z)+\rho^2)=0
    \bmod \sigma(z)  ~~.$$
Since every factor of the product $(a(z)+1)(a(z)+\rho)(a(z)+\rho^2)$ has degree
 less than $t$, at least two of them must have a common non-trivial factor with
  $\sigma(z)$, unless $a(z)=1, \rho, \rho^2$. In this latter case, the Cantor-Zassenhaus 
  algorithm considers another random polynomial instead of $c(z)$, and reiterates the procedure until all factors have been found.

Notice that $a(z)\equiv 0$ never occurs, since $c(z)$ has degree less than $\sigma(z)$, so that at least one
 root of $\sigma(z)$, say $\beta$, is not a root of $c(z)$; then substituting $\beta$ in the identity $c(z)^{\ell_{m}}=q(z)\sigma(z)+a(z)$, we get $a(\beta)\neq 0$, therefore $a(z)$ is not identically zero (this holds even if the roots of $\sigma$ were not in the field of the coefficients, as in the original description of the algorithm).

For the case $p>2$, the procedure is similar: we would consider $\ell_{m} = \frac{p^{m}-1}{2}$
 and $\rho= \alpha^{\ell_m}=-1$, where $\alpha$ is a primitive element of $\mathbb F_{p^m}$.
 Here we would compute $a(z) = c(z)^{\ell_{m}}  \bmod \sigma(z)  $ and then factor as soon as
  $a(z)\neq \pm 1$.

Let us consider now the case $s>1$. One option is to look at
 $\mathbb F_{p^{sm}}$, where the polynomial fully splits into linear factors: once a factor $z-\beta$ is found, it can be multiplied with the factors $z-\beta^{p^{mi}}$, with $1\leq i\leq s-1$, to obtain an irreducible factor of degree $s$. 
A second option is the application of the algorithms over $\mathbb F_{p^{m}}$ (\cite{benor}, \cite{cantor}),
 to directly find the irreducible factors of degree $s$ over $\mathbb F_{p^{m}}$.
If $p=2$, the argument follows as above: either $\gcd(c(z),\sigma(z))$ is non-trivial, or
 $\gcd(c(z),\sigma(z))=1$, in which case 
$$ (c(z)^{\ell_{sm}}+1)(c(z)^{\ell_{sm}}+\rho) (c(z)^{\ell_{sm}}+\rho^2) = (a(z)+1)(a(z)+\rho)(a(z)+\rho^2)=0
    \bmod \sigma(z)  ~~.$$
Since every factor of the product $(a(z)+1)(a(z)+\rho)(a(z)+\rho^2)$ has degree
 less than $t$, at least two of them must have a common non-trivial factor with
  $\sigma(z)$ in $\mathbb F_{2^{m}}$, unless $a(z)=1, \rho, \rho^2$. In this latter case, the Cantor-Zassenhaus algorithm considers another random polynomial $c(z)$, and reiterates the procedure until
  all factors have been found. 

For the case $p>2$, the procedure is similar: we would consider $\ell_{sm} = \frac{p^{sm}-1}{2}$
 and compute 
  $a(z) = c(z)^{\ell_{sm}}  \bmod \sigma(z)$ and then factor as soon as
  $a(z)\neq \pm 1$.
 
In the next Section we will present a variant of the Cantor-Zassenhaus algorithm, according to the description given above, and then deal with probabilistic as well as deterministic considerations about its success rate.

\section{An improved algorithm}

We focus first on the case $s=1$ and show that it is enough, and indeed convenient, to choose $c(z)=z$ as initial test polynomial and to choose $c(z)=z+\beta$, for some random
 $\beta\neq 0$, as further test polynomial, and continuing by choosing random $\beta$s different from the previous ones until a factor is found. A similar approach was already present in \cite{rabin} for the case of odd characteristic (cf. also \cite{shallit}). 

We then consider the case $s>1$, where polynomials of degree $1$ or $s$ will be involved as test polynomials in order to obtain bounds on the number of attempts to find a factor.
\subsection{Case $s=1$}
Suppose $\sigma(z)$ is over $\mathbb F_{2^m}$ and $    z^{\ell_{m}}=\rho^i  \bmod \sigma(z)  $, $i\in\{0,1,2\}$. Now, any element in $\mathbb F_{2^m}^*$ can be written as $\alpha^{k+3n}$, with $k\in\{0,1,2\}$: we define $\mathcal A_0=\{\alpha^{3i}:~~i=0,\ldots , \ell_{m}-1 \}$, that is the subgroup of
 the elements of $\mathbb F_{2^m}^*$ that are cubic powers, and let $\mathcal A_1=\alpha \mathcal A_0$
 and $\mathcal A_2=\alpha^2 \mathcal A_0$ be the two cosets that complete the coset 
 partition of $\mathbb F_{2^m}^*$.  If we substitute $\alpha^{k+3n}$ for any root $z_i$ of $\sigma(z)$ in $z^{\ell_{m}}-\rho^i= Q(z)\sigma(z)$, we obtain $\rho^k-\rho^i=0$, which implies $k=i$. This means that if $    z^{\ell_{m}}=\rho^i  \bmod \sigma(z)  $, then all the roots of $\sigma(z)$ are of the form $\alpha^{i+3n}$, that is they belong to the same coset.
When this situation occurs, we consider another test polynomial $c(z)=z+\beta$, which is equivalent to testing $c(z)=z$ for the polynomial $\varsigma(z)$ whose set of roots is $\{z_i+\beta\}$.
The test succeeds as soon as we find a $\beta$ such that the roots $z_i+\beta$ do not all belong to the same coset. 

The next step is to determine an upper bound to the number of attempts needed in the worst case
 scenario, or on average, until a factor is found. 

Let us first consider the simple case $t=2$:
suppose that $z_1$ and $z_2$ belong to the same coset; then we look for a $\beta$ such that $z_1+\beta$ and $z_2+\beta$ are in different cosets. 
For the worst case scenario, we need to know how many pairs $(z_1+\beta,z_2+\beta)$ have both elements in the same coset. This is equivalent to knowing the number of ways in which $z_1-z_2=z_1+\beta-(z_2+\beta)$
 can be written as the sum of two elements in the same coset. This number is actually $\frac{2^{m}-1}{3}-1$, as can be deduced from \cite[Theorem 1]{winter} specialized with $i=0$ and $\chi$ the cubic character. So at most with $\frac{2^{m}-1}{3}$ attempts we can factor a polynomial of degree $2$. Clearly at each test we can factor with a probability of $\frac{2}{3}$, so that the expected number of attempts is $1.5$.

If $\sigma(z)$ is a polynomial over $\mathbb F_{p^m}$, $p>2$, then the maximum number of attempts is $\frac{p^{m}-1}{2}$, by similar reasoning: we again use some additive properties of residues (\cite{elia2, monico, raymond, winter}). At each test we can factor with a probability of $\frac{1}{2}$, so that the expected number of attempts is $2$.

The remainder of this paper will be devoted to establishing both probabilistic estimates and deterministic bounds on the number of attempts needed to successfully factor, for a generic $t$.
 A first deterministic, though very loose, bound is the following:

\begin{proposition}
   \label{prop1}
The maximum number of attempts needed to find a factor is upper bounded by $\ell_m$ (that is $\frac{2^m-1}{3}$ or $\frac{p^m-1}{2}$ for $p=2$ or $p$ odd, respectively). In particular, in the
 Cantor-Zassenhaus algorithm it is sufficient to consider only linear polynomials as test polynomials $c(z)$.
\end{proposition}

\noindent
{\sc Proof}.
In characteristic $2$, if a root $z_i $ belongs to a given known coset, we can test all the $\ell_m$ elements of that coset, until we obtain $z_i$ itself: $z_i+z_i$ adds to $0$, which does not belong to any coset. Thus we will succeed with at most $\ell_m$ attempts. In characteristic $p$ greater than $2$, it is sufficient to add all the elements of the coset multiplied by $p-1$. 

That it is enough to consider all the $p^m$ monic linear polynomials is anyway clear since computing 
$\gcd \{z-\beta, \sigma(z) \}$ for all $\beta$ in $\mathbb F_{p^m}$ would be enough to find all the factors.
\QED

\paragraph{Remark 1}
The above argument implies that, if the first attempt fails, we know which coset the roots belong to, and can
 restrict our choice of $\beta$ to that coset.
%

\paragraph{Remark 2}
Alternatively, the upper bounds of the proposition follow from the above remarks about $t=2$: clearly, if $t$ is bigger than $2$, then a degree-$2$ polynomial is anyway a factor of the $t$-degree polynomial, so that the maximum number of attempts cannot exceed the number needed to factor this degree-$2$ polynomial.

\paragraph{Remark 3}
In the original version of the Cantor-Zassenhaus algorithm, $\gcd(a(z),\sigma(z))$ is computed
 when searching for a factor of $\sigma(z)$, corresponding to the case when $\gcd(c(z),\sigma(z))$ is non-trivial.
Our version of the algorithm avoids this computation, since it is sufficient to evaluate $\sigma(z)$ in $\beta$ with any efficient polynomial evaluation algorithm; this can be done before exponentiating to the power $\ell_m$.

\paragraph{Remark 4}
If $q$ is a prime factor of $p^m-1$, then we may consider the exponent $\ell_m=\frac{p^m-1}{q}$: in this case the probability of success is $\frac{q-1}{q}$ and the corresponding expected number of attempts is $\frac{q}{q-1}$, which is close to $1$ already for small primes like $5$ or $7$; the drawback is that, if $q$ is large, in the worst case we must check $q$ greatest common divisors, namely $\gcd(a(z)+\zeta_q^j,\sigma(z))$, for $0\leq j\leq q-1$, where $\zeta_q$ is a $q$-th primitive complex root of unity.  

\subsection{Case $s>1$}
If $s>1$, either we look for linear factors in $\mathbb F_{p^{ms}}$, and the analysis is the same as in the case $s=1$, or we choose the direct method, as explained in the previous section. In this case, by a similar argument as above, the algorithm succeeds as soon as $c(z_i)$, $z_i$ being the roots of $\sigma(z)$, are not all in the same coset. This is equivalent to ask that non conjugate roots are not all in the same coset, as 
$$
c(z_i^{p^m})^{\ell_{sm}}=((c(z_i))^{p^m})^{\ell_{sm}}=((c(z_i))^{\ell_{sm}})^{p^m}=(c(z_i))^{\ell_{sm}}
$$
by the properties of the Frobenius automorphism.

Let us see this more precisely, describing in detail the case $p=2$, 
while a similar argument applies in the case of odd primes. 
Let $\sigma(z)$ be, as above, a polynomial of degree $t$ over $\mathbb F_{2^{m}}$, which is a product of $d$
 irreducible polynomials $\sigma_i(z)$ of degree $s$ over the same field $\mathbb F_{2^{m}}$, 
 where it is not restrictive to assume even $m$. 
%
%
According to 
Cantor-Zassenhaus algorithm, a polynomial $c(z)$ over $\mathbb F_{2^m}$,
 relatively prime with $\sigma(z)$, separates $\sigma(z)$ into two polynomials of smaller degree if
 $a(z)=c(z)^{\ell_{sm}}  \bmod \sigma(z)$ is different from $1,\rho,\rho^2$: at 
 least two factors $\sigma_i(z)$ are in two distinct greatest common divisors between $\sigma(z)$
 and $a(z)+1$, $a(z)+\rho$, and $a(z)+\rho^2$, respectively. 
 
\begin{lemma}
  \label{lemma1}
With the above hypotheses and definitions, a polynomial $c(z)$ over $\mathbb F_{2^m}$ separates
 $\sigma(z)$ into two polynomials one containing the factor $\sigma_1(z)$, and a second one containing
 the factor $\sigma_2(z)$ if and only if 
 $c(z)^{\ell_{sm}} \bmod \sigma_1(z)\neq c(z)^{\ell_{sm}} \bmod \sigma_2(z)$. Equivalently, 
 $\sigma_1(z)$ and $\sigma_2(z)$ are separated if and only if $c(z_1)$ and $c(z_2)$ belong to different
 cosets $\mathcal A'_h$ of $\mathbb F_{2^{sm}}^*$, where $z_1$ and $z_2$ are roots of $\sigma_1(z)$ and $\sigma_2(z)$, respectively.
\end{lemma}

\noindent
{\sc Proof}. The polynomial $\sigma(z)$ can be written as a product of three polynomials, i.e.
 $\sigma_1(z)$, $\sigma_2(z)$, and $\sigma_r(z)$ which collects the remaining factors, thus
 $a(z)$ can be decomposed, using the Chinese Remainder Theorem (CRT), as
$$ a(z) = a_1(z) \psi_1(z) + a_2(z) \psi_2(z) + a_r(z) \psi_r(z)   \bmod \sigma(z)  ~~,~~
      \psi_1(z) + \psi_2(z) + \psi_r(z) =1  ~~, $$ 
where $a_1(z) = c(z)^{\ell_{sm}}  \bmod \sigma_1(z)$, $a_2(z) = c(z)^{\ell_{sm}}  \bmod \sigma_2(z)$, 
 and $a_r(z) = c(z)^{\ell_{sm}}  \bmod \sigma_r(z)$.

 If $a(z)=1,\rho, \rho^2$, the uniqueness of the CRT decompositions implies that $a_1(z) = a_2(z)=a_r(z)$. 
 
 If $a(z)\neq 1,\rho, \rho^2$, then $c(z)$ separates $\sigma(z)$
     into two polynomials of smaller degree, and we distinguish two cases:
\begin{itemize}
  \item[1)] $a_1(z) \neq a_2(z)$: 
the polynomials $\sigma_1(z)$ and $\sigma_2(z)$ are in
     different factors because, if both of them were in the same factor, they would both divide the same polynomial
     $a(z)+\rho^h$, thus $a_i(z)=a(z)=\rho^h$ modulo $\sigma_i(z)$, $i=1,2$, contrary to the assumption.
  \item[2)] $a_1(z) = a_2(z)$: 
$\sigma_1(z)$ and
     $\sigma_2(z)$ are in the same factor; in fact, suppose they are not, then 
     $a_1(z) = a(z) =\rho^{h_1} \bmod \sigma_1(z) \neq a_2(z)=a(z) =\rho^{h_2}\bmod \sigma_2(z)$, yielding a contradiction.
\end{itemize}

\noindent
Also, since $ a(z)=c(z)^{\ell_{sm}} \bmod \sigma(z)$ and $a(z) = a_i(z)  = \rho^{h_i} \bmod \sigma_i(z)$, we have that
$c(z_i)^{\ell_{sm}} = \rho^{h_i}$, $i=1,2$, which means that $c(z_i) \in \mathcal A'_{h_i}$,
hence it follows from the first part of the lemma that $c(z)$ separates $\sigma_1(z)$ and $\sigma_2(z)$
 if and only if $c(z_1) \neq c(z_2)$.
\QED

Now, as in the case $s=1$, we are interested in upper bounds for the number of attempts and we can limit the choice of $c(z)$, according to our convenience.
%
%
%
For example, if we know at least one
 primitive polynomial $m(z)$ of degree $s$, we can choose the polynomials $c(z)$ within the set of monic
 irreducible polynomials of degree $s$, so that we get directly $\frac{p^{ms}}{s}$ as an upper bound.
If we do not have any primitive polynomial of degree $s$, that is no means to get and draw from the pool of irreducible polynomials of degree $s$, then we can choose the polynomials $c(z)$ within the larger set of monic polynomials of degree $s$, and we have the looser bound $p^{ms}$.
Somehow surprisingly, we show next that usually it is actually sufficient to consider again linear polynomials.


\noindent

\noindent 
Let $\chi'_3(x)$ be a non-trivial cubic character over $\mathbb F_{2^{sm}}$, namely $\chi'_3$ is a mapping from $\mathbb F_{2^{sm}}^*$ into the complex numbers defined as
$$  \chi'_3(\alpha^h \theta) =  \zeta_3^h ~~~~\theta\in\mathcal A'_0,~~~ h=0,1,2~~,  $$ 
$\alpha$ being a primitive element of $\mathbb F_{2^{sm}}^*$, $\zeta_3$ a primitive complex cubic root of unity, and $\mathcal A'_0$ the coset of cubes in $\mathbb F_{2^{sm}}^*$.
Moreover, we set $\chi'_3(0)=0$ by definition.

If $z_1$ and $z_2$ are roots of two distinct irreducible polynomials of degree $s$, we denote with $N_2^{(m)}(z_1,z_2)$
 the number of monic polynomials $c(z) = z+\beta$ with $\beta\in\mathbb F_{2^m}$ such that 
$ \chi'_3(c(z_1)) =  \chi'_3(c(z_2))$. 
 
%
 
\begin{proposition}
 \label{prop2}
The maximum number $N_A$ of attempts needed to find an irreducible factor of degree $s$, using monic linear polynomials as test polynomials, is upper bounded by
 $\frac{2^{m}}{3}  (1+ \frac{4s-2}{\sqrt{2^m}}+\frac{1}{2^{m}} )$ if $p=2$,
or by $\frac{p^{m}}{2}  (1+ \frac{2s-1}{\sqrt{p^m}})$ if $p$ is odd. In particular linear polynomials are sufficient to find a factor if $\frac{4s-2}{\sqrt{2^m}} <2$ or $ \frac{2s-1}{\sqrt{p^m}} < 1$, respectively.
\end{proposition}

\noindent
{\sc Proof}. 
 In the case of characteristic $2$, $N_A$ is upper bounded by the maximum of $N_2^{(m)}(z_1,z_2)+1$ taken over all distinct pairs of roots $z_1$ and $z_2$ of distinct irreducible polynomials of degree $s$. Thus an upper bound for $N_2^{(m)}(z_1,z_2)$ independent of $z_1$ and $z_2$ is also an upper bound for $N_A-1$. 

 Consider the indicator function
$$I_{\mathcal A'_h}(c(z_i)) = \frac{1+\bar \zeta_3^h \chi'_3(c(z_i))+\zeta_3^h \bar \chi'_3(c(z_i))}{3}
   ~~~~i=1,2 ~~,     $$
which is $1$ if the cubic character of $c(z_i)$ is $\zeta_3^h$, and is $0$ otherwise, if we suppose $c(z)$ relatively prime with $\sigma(z)$.

 Therefore, for a given $c(z)$ we have a coincidence whenever the product 
 $ I_{\mathcal A'_h}(c(z_1)) I_{\mathcal A'_h}(c(z_2)) $ is $1$.
Thus, 
$$  \sum_{h=0}^2 I_{\mathcal A'_h}(c(z_1)) I_{\mathcal A'_h}(c(z_2)) =   \frac{1}{3} \left( 1+ \chi'_3(c(z_1)) \bar \chi'_3(c(z_2)) + \bar \chi'_3(c(z_1)) \chi'_3(c(z_2)) \right)  ~~ $$
is the coincidence indicator for a fixed polynomial $c(z)$.
Summing over all monic linear polynomials $z+\beta$ over $\mathbb F_{2^{m}}$, we get the total number
 $N_2^{(m)}(z_1,z_2)$ of coincidences 
$$ N_2^{(m)}(z_1,z_2)= \frac{1}{3}  \sum_{\beta \in \mathbb F_{2^{m}}} \left( 1+ \chi'_3(z_1+\beta) \bar \chi'_3(z_2+\beta) + \bar \chi'_3(z_1+\beta) \chi'_3(z_2+\beta) \right)-\frac{2}{3}  ~~, $$
where $-\frac{2}{3}$ comes from excluding those polynomials $z+\beta$ having $z_1$ or $z_2$ as root.
We split the summation in three summations, the first summation is simply $2^{m}$, and the second and third summations are complex conjugated, thus it is enough to evaluate only the summation
$$ C= \sum_{\beta \in \mathbb F_{2^{m}}} \chi'_3(z_1+\beta) \bar \chi'_3(z_2+\beta)   ~~.  $$
This summation is hard to evaluate in closed form, thus we content ourselves with a bound. Namely, as  $\chi'_3$ can be considered as the lifted character of a nontrivial character $\chi_3$ over $\mathbb F_{2^{m}}$ \cite{Jungnickel}, we can write
$$
C=\sum_{\beta \in \mathbb F_{2^{m}}} \chi_3(N_{F_{2^{ms}}/F_{2^m}}(z_1+\beta))\bar\chi_3(N_{F_{2^{ms}}/F_{2^m}}(z_2+\beta)),
$$
where $N_{F_{2^{ms}}/F_{2^m}}(x)\doteq x\cdot x^{2^m}\cdots x^{2^{m(s-1)}}$ is the relative norm of $x$.

Since $N_{F_{2^{ms}}/F_{2^m}}(z_i+\beta)$, $i=1,2$, are polynomials of degree $s$ in $\beta$, and $\bar\chi_3=\chi_3^2$, we can then use the Weil bound (\cite[Theorem 2C']{schmidt}; cf. also \cite{wan},\cite[Lemma 2.2]{winterhof}) to obtain 
$$C< (2s-1)2^{m/2}.$$

In conclusion we obtain $N_A$ bounded as
$$  N_A < \frac{2^{m}}{3}  (1+ \frac{4s-2}{\sqrt{2^m}}+\frac{1}{2^{m}} )~~.  $$
The same argument works similarly for $p$ odd, and making the appropriate changes the conclusion is
$$  N_A < \frac{p^{m}}{2}  (1+ \frac{2s-1}{\sqrt{p^m}})~~.  $$
%
\QED

\vspace{5mm}

In the following we analyse the algorithm more in detail both from a probabilistic and a deterministic point of view; in particular we will show that the maximum number of attempts to get a factor is usually very small, so that the algorithm, which is probabilistic in nature, can often be considered deterministic. In order to simplify the subsequent analysis, we will suppose that $s=1$ from now on. 

\section{Probability of factoring}
The Cantor-Zassenhaus algorithm is very efficient in factoring polynomials, but is not deterministic. We can show, however, that the maximum number of attempts, following the modified version above, decreases exponentially with the degree of the polynomial, so that the probability of factoring with one test is close to $1$ when the degree is large enough. 

Making the reasonably assumption that the set of $\{z_i+\beta\}$ for some $\beta$ is made up of elements which belong to each coset $\mathcal A_i$ with probability $1/3$ (or $1/2$ in the case $p>2$), independently of one another, then $3\cdot\frac{1}{3^t}$ is the probability that they all belong to a common coset of the three cosets (and $2\cdot\frac{1}{2^t}$ in case of the two cosets in $\mathbb F_{p^m}^*$, $p>2$).
Therefore the number of attempts to obtain a factor, in the worst case scenario, is roughly $\frac{2^m}{3^{t-1}}$ and $\frac{p^m}{2^{t-1}}$ respectively. And the expected number is $\frac{1}{1-\frac{1}{3^{t-1}}}=1+\frac{1}{3^{t-1}-1}$ or $1+\frac{1}{2^{t-1}-1}$. 

Furthermore, suppose we fail at the first attempt, then we can choose $\beta$ within a certain coset, and the probability of failing at the next $n$ attempts is only $\frac{1}{3^{tn}}$.

Clearly, once a factor is found, the polynomial splits into two parts to which we will re-apply
 the previous computation if we are interested in a complete factorization, untill all linear factors
 are obtained.

\section{Deterministic splitting I: fixed $t$}
If we use the proposed variant of the Cantor-Zassenhaus algorithm, the tightest upper bound to the number of
 attempts necessary to split a polynomial $\sigma(z)$ of degree $t$ over $F_{2^m}$ is equal to 
$$
1+\max_{z_1\neq z_2\neq\cdots\neq z_t} N_2(t),
$$
 where $N_2(t)$ is the number  of
 solutions $\beta$ of a system of $t$ equations in $\mathbb F_{2^m}$ of the form
\begin{equation}
   \label{sys1}
  \left\{ \begin{array}{l}
        \alpha^i z_1^3 +\beta= \alpha^k y_1^3 \\
       \alpha^i  z_2^3 +\beta= \alpha^k y_2^3 \\
 ~~~~~~\vdots    \\
        \alpha^i z_t^3 +\beta= \alpha^k y_t^3 \\ 
    \end{array}  \right.
\end{equation}
where $\alpha^i z_1^3, \alpha^i  z_2^3, \cdots, \alpha^i z_t^3$ are given and distinct (i.e. they are the roots of $\sigma(z)$), whereas the $y_i$s must be chosen in the field to satisfy the system, and the three values $\{0,1,2\}$ for $k$ and $i$ are all considered.
 However, we may assume $i=0$, since dividing each equation by $\alpha^i$, and setting
  $\beta'=\beta \alpha^{-i}$ and $k'=k-i \bmod 3$, we see that the number of solutions 
  of the system is independent of $i$.
If the system is unsolvable, then the number of attempts is $1$.

To evaluate $N_2(t)$, we define an indicator function of the sets $\mathcal A_u$ using the cubic character,
namely for every $x \neq 0$ 
$$    I_{\mathcal A_j}(x) = \frac{1+\zeta_3^{2j}\chi_3(x)+\zeta_3^{j}\bar \chi_3(x)}{3} = \left\{
     \begin{array}{l}
       1 ~~~~\mbox{if}~~ x \in \mathcal A_j   \\
       0 ~~~~\mbox{otherwise}
     \end{array} \right.
        ~~j=0,1,2~~, $$ 
(where the bar denotes complex conjugation). 
Then, given a $z_i$ we can partition the elements $\beta\neq z_i^3$ in $\mathbb F_{2^m}$      
 into subsets depending on the $k\in\{0,1,2\}$ such that $\chi_3(\beta+z_i^3)=\zeta_3^k$. 
Therefore, a solution of (\ref{sys1}) for a fixed $k$ and $i=0$ is singled out by the
 product
$$  \prod_{i=1}^t I_{\mathcal A_k}(\beta+z_i^3)=\frac{1}{3^t} [1+\sum_{i=1}^t \sigma_i^{(k)} ]   ~~, $$
where each $\sigma_i^{(k)}$ is a homogeneous sum of monomials which are products of $i$ characters of the form $\chi_3(\beta+z_h^3)$ or $\bar\chi_3(\beta+z_h^3)$.
Thus $N_2(t)$ is
\begin{equation}
   \label{main1}
   N_2(t)= \sum_{\stackrel{\beta \in \mathbb F_{2^m}}{\beta\not \in\{z_i^3\}}} 
        \left[ \prod_{i=1}^t I_{\mathcal A_0}(\beta+z_i^3) + \prod_{i=1}^t I_{\mathcal A_1}(\beta+z_i^3)+\prod_{i=1}^t I_{\mathcal A_2}(\beta+z_i^3) \right] ~~.
\end{equation}
The roots $z_i$ in the sum need not be considered, since in any case they are not solutions ($z_i^3+z_i^3=0$ cannot be in the same coset as $z_i^3+z_j^3$ if $i\neq j$).

Similarly, in characteristic greater than $2$, the tightest upper bound to the number of
 attempts necessary to split a polynomial $\sigma(z)$ of degree $t$ is equal to 
$$
1+\max_{z_1\neq z_2\neq\cdots\neq z_t} N_p(t),
$$
 where $N_p(t)$ is the number  of
 solutions $\beta$ of a system of $t$ equations in $\mathbb F_{p^m}$ of the form
\begin{equation}
   \label{sys2p}
  \left\{ \begin{array}{l}
        \alpha^i z_1^2 +\beta= \alpha^k y_1^2 \\
       \alpha^i  z_2^2 +\beta= \alpha^k y_2^2 \\
 ~~~~~~\vdots    \\
        \alpha^i z_t^2 +\beta= \alpha^k y_t^2 \\ 
    \end{array}  \right.
\end{equation}
where $\alpha^i z_1^2, \alpha^i  z_2^2, \cdots, \alpha^i z_t^2$ are given and distinct and the two values $\{0,1\}$ for $k$ and $i$ are considered.
 Again, we may assume $i=0$ and we can define an indicator function of the sets $\mathcal B_u$ using the
  quadratic character, where $\mathcal B_0$ is the set of squares and $\mathcal B_1$ the complementary set in $\mathbb F_{p^m}^*$: namely, let $\chi_2$ be a mapping from $\mathbb F_{p^m}^*$ into the complex numbers defined as
$$  \chi_2(\alpha^h \theta) = (-1)^h ~~~~\theta\in\mathcal B_0,~~~~h=0,1~~.  $$ 
 Again, we set $\chi_2(0)=0$.

The corresponding indicator function is thus
$$    I_{\mathcal B_j}(x) = \frac{1+(-1)^{j}\chi_2(x)}{2} = \left\{
     \begin{array}{l}
       1 ~~~~\mbox{if}~~ x \in \mathcal B_j   \\
       0 ~~~~\mbox{otherwise}
     \end{array} \right.
        ~~j=0,1~~. $$ 
Given a $z_i$ we partition $\mathbb F_{p^m}\setminus\{z_i^2\}$ into subsets depending on the value of $k$, such that $\chi_2(\beta+z_i^2)=(-1)^k$. 
Therefore, a solution of (\ref{sys2p}) for a fixed $k$ is given by the product
$$  \prod_{i=1}^t I_{\mathcal B_k}(\beta+z_i^2)=\frac{1}{2^t} [1+\sum_{i=1}^t \sigma_i^{(k)} ]   ~~, $$
where each $\sigma_i^{(k)}$ is a homogeneous sum of monomials which are product of $i$ characters of the form $\chi_2(\beta+z_h^2)$.
Thus $N_p(t)$ is
\begin{equation}
   \label{main2p}
   N_p(t)= \sum_{\stackrel{\beta \in \mathbb F_{p^m}}{\beta\not \in\{-z_i^2\}}} 
        \left[ \prod_{i=1}^t I_{\mathcal B_0}(\beta+z_i^2) + \prod_{i=1}^t I_{\mathcal B_1}(\beta+z_i^2) \right] ~~.
\end{equation}

The following subsections deal with computations of $N_p(t)$ for small values of $t$, then with general bounds on $N_p(t)$.

\subsection{Computations for small $t$}

In the following computations, we will use some properties of nontrivial characters that we briefly
 mention: 
 $  \sum_{x \in \mathbb F_{q}} \chi(x)=0$; if $\beta\neq 0$, then $\sum_{x \in \mathbb F_{q}} \chi(x) \bar \chi(x+\beta) =-1$ (\cite{schip3,winter}). Moreover,
$$
   \sum_{x \in  \mathbb F_{2^{m}}}  \chi_3(x) \chi_3(x+1) = G_{m}(1,\chi) = -(-2)^{m/2},
$$
 with $G_{m}(1,\chi)$ being the Gauss sum (\cite{schip3}). 

We will start with the case $p=2$. First we compute $N_2(2)$, already found above with another technique, then analogously $N_2(3)$.

\paragraph{$t=2$.}  Setting $x_i=\beta+z_i^3$, we have
$$  \prod_{i=1}^2 I_{\mathcal A_h}(x_i) =   \frac{1}{9} \left( 1+ \sigma_1^{(h)}+\sigma_2^{(h)} \right)  ~~h=0,1,2~~, $$
where
$$ \begin{array}{lcl}
    \sigma_1^{(h)}&=& \zeta_3^{2h}\chi_3(x_1)+ \zeta_3^{h}\bar \chi_3(x_1)+\zeta_3^{2h}\chi_3(x_2)+\zeta_3^{h}\bar \chi_3(x_2)   \\
    \sigma_2^{(h)}&=& \zeta_3^{h}\chi_3(x_1)\chi_3(x_2)+ \chi_3(x_1)\bar \chi_3(x_2)+  \bar \chi_3(x_1)\chi_3(x_2)+\zeta_3^{2h}\bar \chi_3(x_1)\bar \chi_3(x_2) 
    \end{array} ~~  $$ 
Since $\sigma_1^{(0)}+\sigma_1^{(1)}+\sigma_1^{(2)}=0$ and $\sigma_2^{(0)}+\sigma_2^{(1)}+\sigma_2^{(2)}= 3(\chi_3(x_1)\bar \chi_3(x_2)+ 
   \bar \chi_3(x_1)\chi_3(x_2))$, the sum of the three products $\prod_{i=1}^2 I_{\mathcal A_k}(x_i)$ is 
 $ \frac{1}{3} \left( 1+\chi_3(x_1)\bar \chi_3(x_2)+  \bar \chi_3(x_1) \chi_3(x_2) \right) $, and thus 
 the sum over $\beta$ in the whole field $\mathbb F_{2^m}$, with the exclusion of $\beta=z_1^3$ and $\beta=z_2^3$, is
$$  N_2(2)=  \frac{1}{3} \left( 2^m-2+ \sum_{\beta \neq z_1^3,z_2^3}  \left(\chi_3(\beta+z_1^3)\bar \chi_3(\beta+z_2^3)+  \bar \chi_3(\beta+z_1^3)\chi_3(\beta+z_2^3)\right)  \right)  ~~. $$ 
Let  $S$ denote the above summation, then $S$ can be evaluated in closed form: by the
 substitution $\beta= z_1^3 +\eta$, since $\chi_3$ is a nontrivial cubic character, we have 
$$  S=  \sum_{\eta \neq 0,z_1^3+z_2^3}  \left(\chi_3(\eta)\bar \chi_3(\eta+z_1^3+z_2^3)+  \bar \chi_3(\eta)\chi_3(\eta+z_1^3+z_2^3)\right) =-2   ~~, $$ 
as the summation of each of the two parts gives $-1$ ($z_1^3+z_2^3 \neq 0$ by hypothesis).
 In conclusion, we have
$$ N_2(2)=  \frac{1}{3} \left( 2^m-4 \right)  ~~, $$
so that
$$
1+\max_{z_1\neq z_2} N_2(2)=\frac{1}{3} \left( 2^m-1 \right)  ~~.
$$

\paragraph{$t=3$.} In this case 
$$  \prod_{i=1}^3 I_{\mathcal A_h}(\beta+z_i^3) = \frac{1}{27} \left( 1+ \sigma_1^{(h)}+\sigma_2^{(h)}+\sigma_3^{(h)} \right)   ~~h=0,1,2~~, $$
where
$$ \begin{array}{lcl}
    \sigma_1^{(h)}&=& \zeta_3^{2h}\chi_3(x_1)+ \zeta_3^{h}\bar \chi_3(x_1)+\zeta_3^{2h}\chi_3(x_2)+\zeta_3^{h}\bar \chi_3(x_2)+\zeta_3^{2h}\chi_3(x_3)+\zeta_3^{h}\bar \chi_3(x_3)   \\
    \sigma_2^{(h)}&=& \zeta_3^{h}\chi_3(x_1)\chi_3(x_2)+ \chi_3(x_1)\bar \chi_3(x_2)+  \bar \chi_3(x_1)\chi_3(x_2)+\zeta_3^{2h}\bar \chi_3(x_1)\bar \chi_3(x_2)+   \\
                & & \zeta_3^{h}\chi_3(x_2)\chi_3(x_3)+ \chi_3(x_2)\bar \chi_3(x_3)+  \bar \chi_3(x_2)\chi_3(x_3)+\zeta_3^{2h}\bar \chi_3(x_2)\bar \chi_3(x_3) + \\
                & & \zeta_3^{h}\chi_3(x_3)\chi_3(x_1)+ \chi_3(x_3)\bar \chi_3(x_1)+  \bar \chi_3(x_3)\chi_3(x_1)+\zeta_3^{2h}\bar \chi_3(x_3)\bar \chi_3(x_1) + \\
    \sigma_3^{(h)}&=& \chi_3(x_1)\chi_3(x_2)  \chi_3(x_3) +  \bar \chi_3(x_1) \bar \chi_3(x_2)  \bar  \chi_3(x_3)+ \zeta_3^{2h}\bar  \chi_3(x_1)\chi_3(x_2)  \chi_3(x_3)+ \\
                & & \zeta_3^{2h}\chi_3(x_1) \bar \chi_3(x_2)  \chi_3(x_3) + \zeta_3^{2h} \chi_3(x_1) \chi_3(x_2) \bar \chi_3(x_3)+ \zeta_3^{h} \bar  \chi_3(x_1) \bar \chi_3(x_2) \chi_3(x_3)+  \\
                & & \zeta_3^{h} \chi_3(x_1) \bar \chi_3(x_2) \bar \chi_3(x_3) + \zeta_3^{h} \bar \chi_3(x_1) \chi_3(x_2) \bar \chi_3(x_3)  
    \end{array} ~~ $$ 
We thus have
$$  \begin{array}{lcl}
      \sigma_1^{0}+\sigma_1^{1}+\sigma_1^{2}&=&0   \\
      \sigma_2^{0}+\sigma_2^{1}+\sigma_2^{2}&=& 3(\chi_3(x_1)\bar \chi_3(x_2)+  \bar \chi_3(x_1)\chi_3(x_2)+\chi_3(x_2)\bar \chi_3(x_3)+\bar \chi_3(x_2)\chi_3(x_3)+ \\
      && ~~ ~~ \chi_3(x_3)\bar \chi_3(x_1)+  \bar \chi_3(x_3)\chi_3(x_1) ) \\ 
      \sigma_3^{0}+\sigma_3^{1}+\sigma_3^{2}&=& 3(\chi_3(x_1)\chi_3(x_2) \chi_3(x_3) + \bar \chi_3(x_1) \bar \chi_3(x_2)  \bar  \chi_3(x_3)) 
      \end{array}  ~~    $$
In the summation over $\beta$ of the sum of the three products, the values of $\beta=z_1^3,z_2^3, z_3^3$
 should be excluded. Thus we must compute 
$$  N_2(3)=  \frac{1}{9} \left( 2^m-3+ \frac{1}{3}\sum_{\beta \neq z_1^3,z_2^3,z_3^3}  \left[(\sigma_2^{0}+\sigma_2^{1}+\sigma_2^{2})+
        (\sigma_3^{0}+\sigma_3^{1}+\sigma_3^{2}) \right]  \right)  ~~. $$ 
Therefore, two types of summations must be evaluated, namely
$$ S_2= \sum_{\beta \neq z_1^3,z_2^3,z_3^3} \chi_3(\beta+	z_1^3)\bar \chi_3(\beta+z_2^3)~~~~\mbox{and}~~~~  S_3= \sum_{\beta \neq z_1^3,z_2^3,z_3^3} \chi_3(\beta+	z_1^3) \chi_3(\beta+z_2^3) \chi_3(\beta+z_2^3)  ~~, $$
 the remaining ones being obtained by symmetry or complex conjugation. Considering $S_2$, and defining for short $y_1=z_2^3+z_3^3$, $y_2=z_1^3+z_3^3$,
  and $y_3=z_2^3+z_1^3$, we have
$$ S_2 = -  \chi_3(y_2)\bar \chi_3(y_1) + \sum_{\beta \neq z_1^3,z_2^3}  \chi_3(\beta+	z_1^3)\bar \chi_3(\beta+z_2^3) =
         - \chi_3(y_2)\bar \chi_3(y_1) + \sum_{x \neq 0,y_3}  \chi_3(x)\bar \chi_3(x+y_3) ~~, $$
thus $S_2=- \chi_3(y_2)\bar \chi_3(y_1)-1$. Considering $S_3$ we have
$$ S_3 = \sum_{\beta \neq z_1^3,z_2^3,z_3^3}  \chi_3(\beta+	z_1^3) \chi_3(\beta+z_2^3) \chi_3(\beta+z_3^3) = 
         \sum_{x \neq 0,y_2,y_3}  \chi_3(x) \chi_3(x+y_3) \chi_3(x+y_2)    $$
thus, with the change of variable  $x=1/z$, since the character is cubic we obtain
$$ S_3 = \sum_{z \neq 0,1/y_2,1/y_3}  \chi_3(1+zy_3) \chi_3(1+zy_2) = \sum_{X \neq 1,0,1+y_3/y_2}  \chi_3(X) \chi_3(X\frac{y_2}{y_3}+1+\frac{y_2}{y_3})   $$
$$  \begin{array}{lcl}
S_3&=& \sty \chi_3(y_2)\bar \chi_3(y_3) \sum_{X \neq 1,0,1+y_3/y_2}  \chi_3(X) \chi_3(X+1+\frac{y_3}{y_2}) \\
   && \\
   &=& \sty -1 + \chi_3(y_2) \bar\chi_3(y_3) \sum_{X \neq 0,1+y_3/y_2}  \chi_3(X) \chi_3(X+1+\frac{y_3}{y_2}) \\
   && \\
   &=& \sty -1 +\bar \chi_3(y_2) \bar \chi_3(y_3) \bar \chi_3(y_1) \sum_{x \in  \mathbb F_{2^m}}  \chi_3(x) \chi_3(x+1) ~~. 
\end{array}      
$$
%
In conclusion, we obtain
$$  \begin{array}{lcl} 
     N_2(3) &=&  \frac{1}{9} \left[ 2^m-11-(-2)^{\frac{m}{2}}[\chi_3(y_1y_2y_3)+\bar \chi_3(y_1y_2y_3)]- \left(\chi_3(y_1y_2^2)+\chi_3(y_1^2y_2)+
          \right. \right. \\
     &&  \left. \left. ~~~~ \chi_3(y_2y_3^2)+\chi_3(y_2^2y_3)+ \chi_3(y_3y_1^2)+\chi_3(y_3^2y_1)   \right)  \right] \\
             \end{array} ~~. $$ 
Note that, if $z_1=0$ (which corresponds to choosing $\beta$ in one particular coset), then $y_2$ and $y_3$ are cubes, and the number of solutions is 
$$  N_2(3) =  \frac{1}{9} \left( 2^m-13-[(-2)^{\frac{m}{2}}+2][\chi_3(y_1)+\bar \chi_3(y_1)] \right) ~~. $$  
      
 Finally we focus our interest on the maximum over the $z_i$ and obtain

  $$   1+\max_{z_1\neq z_2\neq z_3} N_2(3)= \left\{
     \begin{array}{l}
       \frac{1}{9}(2^m+2^{m/2}-2) ~~~~\mbox{for}~~ m/2 \mbox{ even}   \\
       \frac{1}{9}(2^m+2^{m/2+1}+1) ~~~~\mbox{for}~~ m/2 \mbox{ odd}
     \end{array} \right.
        ~~. $$      
\vspace{5mm}

Let us deal now with the case $p>2$:

\paragraph{$t=2$.}  In this case, we have
$$  \prod_{i=1}^2 I_{\mathcal B_h}(\beta+z_i^2) =   \frac{1}{4} \left( 1+ \sigma_1^{(h)}+\sigma_2^{(h)} \right)  ~~h=0,1~~, $$
where $\sigma_1^{(h)}= (-1)^{h}\chi_2(x_1)+ (-1)^{h} \chi_2(x_2)$, and $ \sigma_2^{(h)} = \chi_2(x_1)\chi_2(x_2)$. 

Since $\sigma_1^{(0)}+\sigma_1^{(1)}=0$ and $\sigma_2^{(0)}+\sigma_2^{(1)}= 2(\chi_2(x_1) \chi_2(x_2))$, 
  the sum over $\beta$
 in the whole field $\mathbb F_{p^m}$ with the exclusion of $\beta=-z_1^2$ and $\beta=-z_2^2$ is
$$  N_p(2)=  \frac{1}{2} \left( p^m-2 + \sum_{\beta \neq -z_1^2,-z_2^2}  \left(\chi_2(\beta+z_1^2)\chi_2(\beta+z_2^2)\right)  \right)  ~~. $$ 
Let  $S$ denote the above summation: we evaluate it in closed form by substituting
 $\beta=\eta-z_1^2$; since $\chi_2$ is a nontrivial quadratic character, we have 
$$  S=  \sum_{\eta \neq 0,z_1^2-z_2^2}  \left(\chi_2(\eta) \chi_2(\eta+z_2^2-z_1^2)\right) =-1   ~~, $$ 
  the summation being independent of the term $z_2^2-z_1^2$, which is non-zero by hypothesis.
  In conclusion we have
$$ N_p(2)=  \frac{1}{2} \left( p^m-3 \right)  ~~, $$
so that
$$
1+\max_{z_1\neq z_2} N_p(2)=\frac{1}{2} \left( p^m-1 \right)  ~~.
$$

\paragraph{$t=3$.} In this case 
$$  \prod_{i=1}^3 I_{\mathcal B_h}(\beta+z_i^2) = \frac{1}{8} \left( 1+ \sigma_1^{(h)}+\sigma_2^{(h)}+\sigma_3^{(h)} \right)   ~~h=0,1~~, $$
where $\sigma_1^{(h)}=  (-1)^{h} \chi_2(x_1)+(-1)^{h}\chi_2(x_2)+(-1)^{h}\chi_2(x_3)$, 
      $\sigma_2^{(h)}=  \chi_2(x_1)\chi_2(x_2)+ \chi_2(x_1)\chi_2(x_3)+ \chi_2(x_2) \chi_2(x_3) $, and
      $\sigma_3^{(h)}=  (-1)^{h} \chi_2(x_1)\chi_2(x_2) \chi_2(x_3)$.

Since $\sigma_1^{0}+\sigma_1^{1}=0 $, $\sigma_2^{0}+\sigma_2^{1}= 2(\chi_2(x_1) \chi_2(x_2)+ \chi_2(x_1)\chi_2(x_3)+\chi_2(x_2) \chi_2(x_3) )$, and $\sigma_3^{0}+\sigma_3^{1}= 0$,
 the summation over $\beta$ of the sum of the two products, where the values of $\beta$ equal to $-z_1^2,-z_2^2$, and $-z_3^2$ are excluded, becomes 
$$  N_p(3)=  \frac{1}{4} \left( p^m-3+ \sum_{\beta \neq -z_1^2,-z_2^2,-z_3^2}  \left[ \chi_2(x_1) \chi_2(x_2)+ \chi_2(x_1)\chi_2(x_3)+\chi_2(x_2) \chi_2(x_3)  \right]  \right)  ~~. $$ 
We thus need to evaluate only one type of summation, namely
$$ S_2= \sum_{\beta \neq -z_1^2,-z_2^2,-z_3^2} \chi_2(\beta+	z_1^2) \chi_2(\beta+z_2^2) = 
   \sum_{\stackrel{\eta \neq 0}{z_1^2-z_2^2,z_1^2-z_3^2}} \chi_2(\eta) \chi_2(\eta+z_2^2-z_1^2) =-1- \chi_2(z_1^2-z_3^2) \chi_2(z_2^2-z_3^2)~, $$
 the remainder being obtained by symmetry. In conclusion, we obtain
$$  N_p(3) =  \frac{1}{4} \left[ p^m-6  - (\chi_2(z_1^2-z_3^2) \chi_2(z_2^2-z_3^2)+\chi_2(z_1^2-z_2^2) \chi_2(z_3^2-z_2^2)+
 \chi_2(z_3^2-z_1^2) \chi_2(z_2^2-z_1^2)) \right]  ~~. $$ 

And , if we consider the maximum, we have

$$
1+\max_{z_1\neq z_2\neq z_3} N_p(3)==\left\{\begin{array}{l}
\frac{1}{4}(p^m-1)~~p=4k+1\\
\frac{1}{4}(p^m+1)~~p=4k+3,~~m\ \mbox{odd}\\
\frac{1}{4}(p^m-1)~~p=4k+3,~~m\ \mbox{even}
\end{array}\right.
$$

\subsection{Bounds}
As the number of equations in system \ref{sys1} or \ref{sys2p} becomes larger, exact computations
 become less meaningful for our purpose, as it would then be necessary to think about estimates
 and bounds on rather cumbersome expressions. We will thus shift our interest to a general upper
  bound for the function $N_p(r)$; we will first deal with the case $p=2$, then
  the case $p>2$. 
 
Consider equation (\ref{main1}) written as
\begin{equation}
   N_2(r)= \frac{1}{3^r} \sum_{\stackrel{\beta \in \mathbb F_{2^m}}{\beta\not \in\{z_i^3\}}} 
         \left[ \mathfrak P_0 +  \mathfrak P_1+ \mathfrak P_2 \right] ~~,
\end{equation}
where
$$   \mathfrak P_k = 3^r \prod_{i=1}^r I_{\mathcal A_k}(x_i)=  1+\sigma_1^{(k)}+ \sigma_2^{(k)} + \cdots + \sigma_r^{(k)}  ~~~~~~k=0,1,2 ~~,$$
$x_i$ being $\beta+z_i^3$, and each $\sigma_j^{(k)}$ is a sum of monomials  which are products of the same number $j$ of distinct variables
 (characters) $\chi_3(x_i)$ or $\bar\chi_3(x_i)$, possibly times $\zeta_3$ or $\zeta_3^2$. 
In particular the number of addends in $\sigma_j^{(k)}$ is $2^j \binomiale{r}{j}$. 

Define $\sigma_j=\sigma_j^{(0)}+\sigma_j^{(1)}+\sigma_j^{(2)}$ for every $j=1, \ldots, r$; then $\sigma_j$
 contains fewer addends than any $\sigma_j^{(k)}$, since all monomials multiplied by either $\zeta_3$ or $\zeta_3^2$ are canceled out with monomials multiplied by $1$, and the surviving monomials are multiplied by
  $3$ (see also the examples above). 
In particular, $\sigma_1$ is zero; $\sigma_2$ is a sum of monomials of the form $\chi_3(x_i) \bar \chi_3(x_l)$
  ($i,l$ distinct),
  whose total number is $2\binomiale{r}{2}$; $\sigma_3$ is a sum of monomials of the form $\chi_3(x_i) \chi_3(x_l)\chi_3(x_m)$ ($i,l,m$ all distinct), whose total number is $2\binomiale{r}{3}$; and 
  $\sigma_4$ is a sum of monomials of the form $\chi_3(x_i)\chi_3(x_l)\bar \chi_3(x_m)\bar \chi_3(x_s)$ ($i,l,m,s$ all distinct), whose  total number is $6\binomiale{r}{4}$.
%
In general, the number of surviving monomials of degree $j$ can be computed by considering that
 each monomial is a product of $n_1$ characters and $n_2$ complex conjugate characters; thus $n_1+n_2=j$.
Supposing that $\chi_3(x_i)$ are multiplied by $\zeta_3$ and  $\bar \chi_3(x_h)$ are multiplied by $\zeta_3^2$, 
 the surviving monomial satisfies the condition $n_1+2 n_2= 0 \bmod 3$. Therefore, the admissible values 
 of $0 \leq n_2 \leq j$ satisfy the condition $n_2= 2j \bmod 3$: if $e= 2j \bmod 3$ and $e\in\{0,1,2\}$,
 the number of surviving monomials is
$  \binomiale{r}{j} a_j, $  
where $a_j=\sum_{h=0}^{\lfloor \frac{j-e}{3} \rfloor} \binomiale{j}{e+3h}$,
%
with $\{a_j\}_{\mathbb Z_{>1}}=2,2,6,10,22,42,86,170,342 \ldots$  
 matching the sequence A078008 in \cite{sloane} with the first two terms disregarded.
We observe now that the product of $j$ characters, whose arguments are distinct linear functions of
 $\beta$, can be interpreted as a single character whose argument is a polynomial $f(\beta)$ with $j$ distinct
 roots: by \cite[Theorem 2C']{schmidt}, each sum of these characters is upper bounded by
 $(j-1) \sqrt{2^m}$, so that 
$$  N_2(r)\leq \frac{1}{3^{r-1}} \left[2^m-r + \sum_{j=2}^r a_j (j-1) \binomiale{r}{j} \sqrt{2^m} \right] ~~~~.  $$
The summation above is evaluated as follows, using the expression 
$a_j= \frac{1}{3} \sum_{h=0}^2 \zeta_3^{-he}(1+\zeta_3^h)^j$ for the sequence $a_j$ as can be found in \cite{benjamin1,benjamin2,gould}:
$$ \sum_{j=2}^r a_j (j-1) \binomiale{r}{j}  = \sum_{j=2}^r \frac{1}{3} \sum_{h=0}^2 \zeta_3^{-he}(1+\zeta_3^h)^j (j-1) \binomiale{r}{j} = 
\frac{1}{3} \sum_{h=0}^2 \sum_{j=2}^r \zeta_3^{-he}(1+\zeta_3^h)^j (j-1) \binomiale{r}{j}.
$$
Now, observing that $e=-j \bmod 3$ and $\zeta_3$ is a cubic root of the unity, we may substitute $\zeta_3^{hj}$ for $\zeta_3^{-he}$ and write $(\zeta_3^h+\zeta_3^{2h})^j$ for $\zeta_3^{hj}(1+\zeta_3^h)^j$  in the last expression, which we then write as 
$$ \frac{1}{3} \sum_{h=0}^2 \sum_{j=0}^r (\zeta_3^h+\zeta_3^{2h})^j (j-1) \binomiale{r}{j}+1=
 1+\frac{1}{3} \sum_{h=0}^2 \left( \sum_{j=0}^r j(\zeta_3^h+\zeta_3^{2h})^j  \binomiale{r}{j}-
  \sum_{j=0}^r (\zeta_3^h+\zeta_3^{2h})^j \binomiale{r}{j} \right) ~~.
$$ 
Using the binomial sum and its derivative, we finally obtain
$$ \sum_{j=2}^r a_j (j-1) \binomiale{r}{j}  =1+  \frac{1}{3} \sum_{h=0}^2 \left( r(\zeta_3^h+\zeta_3^{2h})(1+\zeta_3^h+\zeta_3^{2h})^{r-1} - (1+\zeta_3^h+\zeta_3^{2h})^r \right) ~~,
$$ 
that is
$$ \sum_{j=2}^r a_j (j-1) \binomiale{r}{j}  =1+  \frac{1}{3} [2r3^{r-1}-3^r  ]  ~~, $$
because $(1+\zeta_3^h+\zeta_3^{2h})$ is $3$ when $h=0$ and is $0$ otherwise.
In conclusion 
$$  N_2(r)\leq \frac{1}{3^{r-1}} \left[2^m+ \sqrt{2^m}-r + 3^{r-2} (2r-3) \sqrt{2^m} \right] ~~~~,  $$
where we see that, when  $3^{r-2} (2r-3)\sqrt{2^m}-r+\sqrt{2^m} <<2^m$, roughly $r<<m/2$, then 
 $ N_2(r)\simeq \frac{2^m}{3^{r-1}}$, so that this deterministic bound supports the probabilistic estimate
 discussed above.

\vspace{3mm}

In the case $p > 2$, consider equation (\ref{main2p}) written as
\begin{equation}
   \label{bbq2}
   N_p(r)= \frac{1}{2^r} \sum_{\stackrel{\beta \in \mathbb F_{p^m}}{\beta\not \in\{-z_i^2\}}} 
         \left[ \mathfrak Q_0 +  \mathfrak Q_1 \right] ~~,
\end{equation}
where
$$   \mathfrak Q_k = 2^r \prod_{i=1}^r I_{\mathcal B_k}(x_i)=  1+\sigma_1^{(k)}+ \sigma_2^{(k)} + \cdots + \sigma_r^{(k)}  ~~~~~~k=0,1 ~~,$$
$x_i$ being $\beta+z_i^2$, and each $\sigma_j^{(k)}$ is a sum of monomials  which are products of the same number $j$ of distinct variables (characters) $\chi_2(x_i)$. 
 In particular, only $\sigma_j^{(k)}$s with even subscripts occur, and clearly they are the elementary symmetric functions of $r$ variables; thus the number of addends in $\sigma_j^{(k)}$ is $ \binomiale{r}{j}$. 
 The same argument used to upper bound $N_2(r)$ also applies here, in this case the sum of products of $j$ characters is bounded as $(j-1) \sqrt{p^m}$ by \cite[Theorem 2C']{schmidt}, so that
$$  N_p(r)\leq \frac{1}{2^{r-1}} \left[p^m-r + \sum_{j=2}^r (j-1) \binomiale{r}{j} \sqrt{p^m} \right] ~~~~.  $$ 
which, after some manipulation,  can be written as
$$  N_p(r)\leq \frac{1}{2^{r-1}} \left[p^m-r + [2^{r-1}(r-2)+1]  \sqrt{p^m} \right] ~~~~,  $$
and we see that, when  $[2^{r-1}(r-2)+1]\sqrt{p^m}-r <<p^m$, roughly $r<<\frac{m}{2} \log_2 p$, then $ N_p(r)\simeq \frac{p^m}{2^{r-1}}$ as in our probabilistic estimate.

\section{Deterministic splitting II: fixed $N$}

This section examines the smallest $t$ such that the algorithm succeeds, in at most $1$ or $2$ attempts:
 we will call these $t_0(1)$ and $t_0(2)$, respectively.

Clearly, $t_0(1)=\ell_m+1$, since there are exactly $\ell_m$ elements belonging to a given coset; then, if $t>\ell_m$, the algorithm succeeds at the first attempt.

To evaluate $t_0(2)$, we must examine the number of representations of a $\beta\neq 0$ in the field 
 being the sum of an element in a given coset and an element in another (possibly the same) given coset (see also \cite{monico, elia2, raymond}).
We then consider the maximum $M$, over $\beta\neq 0$ in the field and over all possible pairs of cosets, so that $t_0(2)$ is $1+M$.

For the case of the cubic character, $M$ can be calculated as follows:
$$
M=\max_{i,j,\beta} \sum_{z\neq 0,\beta}\frac{1+\zeta_3^{2j}\chi_3(z)+\zeta_3^{j}\bar \chi_3(z)}{3}~~\frac{1+\zeta_3^{2i}\chi_3(\beta+z)+\zeta_3^{i}\bar \chi_3(\beta+z)}{3}
$$
which is the maximum over $i,j,\beta$ of the following expression:
$$
\frac{1}{9}\left[2^m-2-\chi_3(\beta)(\zeta_3^{2i}+\zeta_3^{2j})-\bar\chi_3(\beta)(\zeta_3^{i}+\zeta_3^{j})-\zeta_3^{2i+j}-\zeta_3^{i+2j}-(-2)^{m/2}(\zeta_3^{2i+2j}\bar\chi_3(\beta)+\zeta_3^{i+j}\chi_3(\beta))\right],
$$
where we have again exploited the relations $  \sum_{x \in \mathbb F_{2^m}} \chi_3(x)=0$, $\sum_{x \in \mathbb F_{2^m}} \chi_3(x) \bar \chi_3(x+\beta) =-1$ and $\sum_{x \in  \mathbb F_{2^{m}}}  \chi_3(x) \chi_3(x+1) = G_{m}(1,\chi_3) = -(-2)^{m/2}$ (\cite{berndt,schip3,winter}). Then we have
$$    M =  \left\{
     \begin{array}{l}
       \frac{1}{9}(2^m+2^{m/2}-2) ~~~~\mbox{for}~~ m/2 \mbox{ even}   \\
       \frac{1}{9}(2^m+2^{m/2+1}+1) ~~~~\mbox{for}~~ m/2 \mbox{ odd}
     \end{array} \right.
        ~~. $$ 

For the case of the quadratic character, we consider similarly 

\vspace{4mm}
\noindent
$$
\begin{array}{lcl}
 \sty
M &=& \sty \max_{i,j,\beta} \sum_{z\neq 0,\beta}\frac{1+(-1)^{j}\chi_2(z)}{2}~\frac{1+(-1)^{i}\chi_2(\beta -z)}{2} \\
 && \\
 &=& \sty \max_{i,j,\beta} \left\{\frac{1}{4}\left(p^m-2-\chi_2(\beta)(-1)^i-\chi_2(\beta)(-1)^j-(-1)^{i+j}\chi_2(-1)\right)\right\},
 \end{array}
$$

 
%

\noindent
therefore 
$$ M=\left\{\begin{array}{l}
\frac{1}{4}(p^m-1)~~p=4k+1\\
\frac{1}{4}(p^m+1)~~p=4k+3,~~m\ \mbox{odd}\\
\frac{1}{4}(p^m-1)~~p=4k+3,~~m\ \mbox{even}
\end{array}\right.
$$

\paragraph{Remark 5}
It is interesting to notice that $M$, which is the maximum $t$ such that it is still possible to fail splitting a polynomial of degree $t$ with two attempts, is equal to the maximum number of attempts to split a polynomial of degree $3$. Similarly, $\ell_m$ is at the same time the maximum $t$ such that it is possible to fail splitting a polynomial of degree $t$ at the first attempt and the maximum number of attempts to split a polynomial of degree $2$. 


\section{Acknowledgments}
We would like to thank Joachim Rosenthal and Elisa Gorla for support and fruitful discussions.

The Research was supported in part by the Swiss National Science
Foundation under grants No. 126948 and 132256.

\end{document}